\documentclass[letterpaper, 10pt, conference]{ieeeconf}
\IEEEoverridecommandlockouts 
\usepackage{amsmath,amssymb,url}
\usepackage{graphicx,subfigure,warmread}
\usepackage[all,import]{xy}

\newcommand{\bracket}[1]{\ensuremath{\left[ #1 \right]}}
\newcommand{\braces}[1]{\ensuremath{\left\{ #1 \right\}}}
\newcommand{\parenth}[1]{\ensuremath{\left( #1 \right)}}

\newcommand{\refeqn}[1]{(\ref{eqn:#1})}
\newcommand{\reffig}[1]{Fig. \ref{fig:#1}}
\newcommand{\tr}[1]{\mbox{tr}\ensuremath{\negthickspace\bracket{#1}}}

\newcommand{\G}{\ensuremath{\mathsf{G}}}
\newcommand{\SO}{\ensuremath{\mathsf{SO(3)}}}
\newcommand{\T}{\ensuremath{\mathsf{T}}}
\renewcommand{\L}{\ensuremath{\mathsf{L}}}
\newcommand{\so}{\ensuremath{\mathfrak{so}(3)}}
\newcommand{\SE}{\ensuremath{\mathsf{SE(3)}}}

\renewcommand{\Re}{\ensuremath{\mathbb{R}}}

\newcommand{\D}{\ensuremath{\mathbf{D}}}
\newcommand{\pair}[1]{\ensuremath{\left\langle #1 \right\rangle}}

\newcommand{\Ad}{\ensuremath{\mathrm{Ad}}}
\newcommand{\ad}{\ensuremath{\mathrm{ad}}}
\newcommand{\g}{\ensuremath{\mathfrak{g}}}

\title{\LARGE \bf
Computational Geometric Optimal Control of\\ Connected Rigid Bodies in a Perfect Fluid}

\author{Taeyoung Lee, Melvin Leok\authorrefmark{1}, and N. Harris McClamroch\authorrefmark{2}%
\thanks{Taeyoung Lee, Mechanical and Aerospace Engineering, Florida Institute of Technology, Melbourne, FL 39201 {\tt taeyoung@fit.edu}}%
\thanks{Melvin Leok, Mathematics, University of California, San Diego, CA 92093 {\tt mleok@math.ucsd.edu}}%
\thanks{N. Harris McClamroch, Aerospace Engineering, University of Michigan, Ann Arbor, MI 48109 {\tt
nhm@umich.edu}}%
\thanks{\textsuperscript{\footnotesize\ensuremath{*}}This research has been supported in part by NSF under grants DMS-0714223, DMS-0726263 and DMS-0747659.}
\thanks{\textsuperscript{\footnotesize\ensuremath{\dagger}}This research has been supported in part by NSF under grant CMS-0555797.}
}

\begin{document}
\allowdisplaybreaks
\maketitle \thispagestyle{empty} \pagestyle{empty}

\begin{abstract}
This paper formulates an optimal control problem for a system of rigid bodies that are connected by ball joints and immersed in an irrotational and incompressible fluid. The rigid bodies can translate and rotate in three-dimensional space, and each joint has three rotational degrees of freedom. We assume that internal control moments are applied at each joint. We present a computational procedure for numerically solving this optimal control problem, based on a geometric numerical integrator referred to as a Lie group variational integrator. This computational approach preserves the Hamiltonian structure of the controlled system and the Lie group configuration manifold of the connected rigid bodies, thereby finding complex optimal maneuvers of connected rigid bodies accurately and efficiently. This is illustrated by numerical computations.
\end{abstract}

\section{Introduction}

Fish locomotion involves a deformable fish body interacting with an unsteady fluid, through which  internal muscular forces on the fish are translated into fish motions~\cite{Sfa.IJoOE1999}.

The planar articulated rigid body model has become popular in engineering, as it describes underwater robotic vehicles that move and steer by changing their shape~\cite{Bar.1996}. Furthermore, if it is assumed that the ambient fluid is incompressible and irrotational, then equations of motion of the articulated rigid body can be derived without explicitly incorporating  fluid variables~\cite{Kan.JoNS2005}. The effect of the fluid is accounted by added inertia terms of the rigid body. This model is known to characterize the qualitative behavior of fish swimming~\cite{Kan.JoNS2005}. 

In~\cite{LeeLeoPACC09}, an analytical model and a geometric numerical integrator for three-dimensional connected rigid bodies immersed in an incompressible and irrotational fluid were developed. The connected rigid bodies can freely translate and rotate in three-dimensional space, and each joint has three rotational degrees of freedom. The geometric numerical integrator presented in~\cite{LeeLeoPACC09} is referred to as a Lie group variational integrator~\cite{Lee.2008}, and it preserves the symplecticity, momentum map, and Lie group configuration manifold of the connected rigid bodies. These properties are important for accurately and efficiently studying complex maneuvers of rigid bodies~\cite{Lee.CMaDA2007}. 

This paper formulates an optimal control problem for connected rigid bodies in a perfect fluid. We assume that internal control moments are applied at each joint. These control moments change the relative attitude between rigid bodies, thereby controlling the shape of a system of connected rigid bodies. By using the nonlinear coupling, referred to as a geometric phase~\cite{Mar92}, between shape changes and group motions, and by using the momentum exchange between the rigid bodies and the ambient fluid, the system of connected rigid bodies can translate and rotate without external forces or moments. 

We present a computational approach to solve optimal control problems based on the Lie group variational integrators developed in~\cite{LeeLeoPACC09}. The optimal control problems are formulated as discrete-time optimal control problems using Lie group variational integrators; these problems can be solved using standard numerical optimization techniques. This is in contrast to conventional approaches where discretization is introduced at the last stage in order to solve the optimality conditions numerically. 

This computational geometric optimal control approach preserves the Hamiltonian structure of the controlled system and Lie group structure of connected rigid bodies in fluid~\cite{LeeLeoCISSIDRWB09}. As it does not introduce artificial numerical dissipation that typically appear in general-purpose numerical integration methods, it improves the overall computational efficiency and accuracy of the optimization process. As it is coordinate-free, this approach avoids singularities and complexities associated with local coordinates.

Optimal shape changes for a planar articulated body to achieve a desired locomotion objective have been studied in~\cite{Kan.PotICoDaC2005,Ros.PotACC2006}. The main contribution of this paper is to develop a computational framework for studying optimal maneuvers of connected rigid bodies in a three-dimensional space. By considering rigid body motions in a Lie group, this approach can find nontrivial optimal maneuvers involving large three-dimensional rotations.

\section{Dynamics of Connected Rigid Bodies Immersed in a Perfect Fluid}\label{sec:CRB}

Consider three connected rigid bodies immersed in a perfect fluid. We assume that these rigid bodies are connected by a ball joint that has three rotational degrees of freedom, and the fluid is incompressible and irrotational. We also assume each body has neutral buoyancy.

\renewcommand{\xyWARMinclude}[1]{\includegraphics[width=0.78\columnwidth]{#1}}
\begin{figure}
$$\begin{xy}
\xyWARMprocessEPS{RBPF3d}{eps}
\xyMarkedImport{}
\xyMarkedMathPoints{1-13}
\end{xy}\vspace*{-0.35cm}
$$

\vspace*{-0.35cm}

\caption{Connected rigid bodies immersed in a perfect fluid}\label{fig:CRB}
\end{figure}
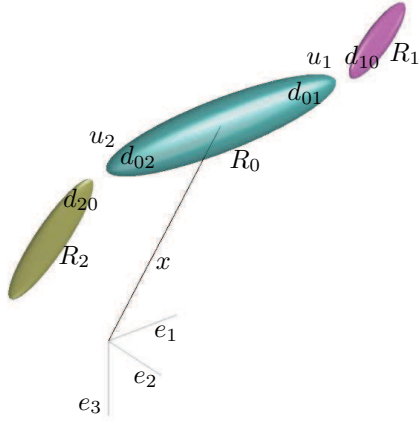

\subsection{Configuration Manifold}

We choose a reference frame and three body-fixed frames. The origin of each body-fixed frame is located at the center of mass of the rigid body. Define
{\allowdisplaybreaks
\begin{center}
\begin{tabular}{lp{6cm}}
$R_i\in\SO$ & Rotation matrix from the $i$-th body-fixed frame to the reference frame\\
$\Omega_i\in\Re^3$ & Angular velocity of the $i$-th body, represented in the $i$-th body-fixed frame\\
$x\in\Re^3$ & Vector from the origin of the reference frame to the center of mass of the $0$-th body, represented in the reference frame\\
$V_i\in\Re^3$ & Velocity of the $i$-th body represented in the $i$-th body fixed frame\\
$d_{ij}\in\Re^3$ & Vector from the center of mass of the $i$-th body to the ball joint connecting the $i$-th body with the $j$-th body, represented in the $i$-th body-fixed frame\\
$m^b_i\in\Re$ & Mass of the $i$-th body\\
$J^b_i\in\Re^{3\times3}$ & Inertia matrix of the $i$-th body\\
$u_i\in\Re^3$ & Internal control moment applied at the $i$-th joint represented in the $0$-th body fixed frame\\
\end{tabular}
\end{center}}
\noindent for $i,j\in\{0,1,2\}$.

A configuration of this system can be described by the location of the center of mass of the central body, and the attitude of each rigid body with respect to the reference frame. The configuration manifold is $\G=\SE\times\SO\times\SO$, where $\SO=\{R\in\Re^{3\times 3}\,|\, R^TR=I, \det{R}=1\}$, and $\SE=\SO\textcircled{s}\Re^3$. 

The attitude kinematic equation is given by
\begin{align*}
    \dot R_i = R_i\hat\Omega_i
\end{align*}
for $i\in\{0,1,2\}$, where the \textit{hat map} $\hat\cdot:\Re^3\rightarrow\so$ is defined by the condition that $\hat x y=x\times y$ for any $x,y\in\Re^3$. The velocity of the central body is $V_0=R_0^T\dot x $. Since the location of the center of mass of the $i$-th rigid body can be written as $x+ R_0d_{0i}-R_id_{i0}$ for $i\in\{1,2\}$ with respect to the reference frame, $V_i$ is given by
\begin{align}
    V_i = R_i^T \dot x - R_i^T R_0 \hat d_{0i}\Omega_0 + \hat d_{i0}\Omega_i.\label{eqn:Vi}
\end{align}

\subsection{Lagrangian}\label{sec:AM}

The total kinetic energy of the connected rigid bodies in a perfect fluid is the sum of the kinetic energy of the rigid bodies and the kinetic energy of the fluid. As the fluid is irrotational, the kinetic energy of the fluid can be expressed without explicitly incorporating fluid variables. The effects of the ambient fluid is accounted for by added inertia terms~\cite{Kan.JoNS2005}. The resulting model is shown to capture the qualitative properties of the interaction between rigid body dynamics and fluid dynamics correctly~\cite{Kan.JoNS2005,Kan.PotICoDaC2005}. 

Let $M^f_{i},J^f_{i}\in\Re^{3\times 3}$ be the added inertia matrices for the translational kinetic energy and the rotational kinetic energy of the $i$-th rigid body. To simplify expressions for the added inertia terms, we assume each body is an ellipsoid. The added inertia matrices depends on the length of principal axes, and their explicit expressions are available in~\cite{Hol.PD1998}.

The total kinetic energy can be written as
\begin{align}
    T & = \sum_{i=0}^2 \frac{1}{2}M_i V_i\cdot V_i + \frac{1}{2} \Omega_i \cdot J_i\Omega_i,
\end{align}
where $M_i = m^b_iI_{3\times 3}  + M^f_{i}$, $J_i  = J^b_i + J^f_{i}$ for $i=\{0,1,2\}$. This can be written in matrix form as
\begin{align}
T = \frac{1}{2} \xi^T \mathbb{I}(R_0,R_1,R_2) \xi,\label{eqn:T}
\end{align}
where $\xi=[\Omega_0;\dot x;\Omega_1;\Omega_2]\in\Re^{12}$ and the matrix $\mathbb{I}(R_0,R_1,R_2)\in\Re^{12\times 12}$~\cite{LeeLeoPACC09}. Since there is no potential field, this corresponds to the Lagrangian of the connected rigid bodies in a perfect fluid.

\subsection{Euler-Lagrange Equations}

Euler-Lagrange equations for a mechanical system that evolves on an arbitrary Lie group are given by
\begin{gather}
\frac{d}{dt}\D_\xi L(g,\xi)-\ad^*_\xi \cdot \D_\xi L(g,\xi) -\T_e^*\L_g\cdot \D_g L(g,\xi)=U,\label{eqn:EL0}\\
\dot g= g\xi,\label{eqn:EL1}
\end{gather}
where $L:\T\G\simeq \G\times\g\rightarrow \Re$ is the Lagrangian of the system~\cite{Lee.2008}. Here $\D_\xi L(g,\xi)\in\g^*$ denotes the derivative of the Lagrangian with respect to $\xi\in\g$, $\ad^*:\g\times\g^*\rightarrow\g^*$ is the $\mathrm{co}$-$\mathrm{adjoint}$ operator, and $\T_e^*\L_g:\T^*\G\rightarrow\g^*$ denotes the cotangent lift of the left translation map, $\L_g:\G\rightarrow\G$~\cite{Mar.BK1999}, restricted to the fiber over the identity $e$.

For the configuration manifold $\G=\SO\times\Re^3\times\SO\times\SO$, we left-trivialize $\T\G$ to yield $\G\times\g$, and we identify its Lie algebra $\g$ with  $\Re^{12}$ by the hat map. Using this result, Euler-Lagrange equations for the uncontrolled dynamics of connected rigid bodies in a perfect fluid have been developed in~\cite{LeeLeoPACC09}. Here, we focus on finding the effect of internal control moments represented by $U\in\g^*$ in \refeqn{EL0}.

The central $0$-th rigid body is subject to the control moment $u_1+u_2$. Let $\rho\in\Re^3$ be the vector  from the center of mass of the $0$-th rigid body to a mass element, and let $dF(\rho)\in\Re^3$ be the force acting on the mass element expressed in the $0$-th body fixed frame. Then, $\int_{\mathcal{B}_0} \rho\times dF(\rho)=u_1+u_2$. The virtual displacement of the mass element due to the rotation of the $0$-th body can be written as $\delta R_0 \rho  = R_0 \hat \eta_0 \rho$ in the reference frame, where the variation of the rotation matrix is represented by $\delta R_0 = \frac{d}{d\epsilon}\big|_{\epsilon=0} R_0 \exp \epsilon \hat\eta_0=R_0\hat\eta_0$ for $\eta_0\in\Re^3$. Then, the virtual work done by the control moment on the $0$-th rigid body is given by
\begin{align*}
\int_{\mathcal{B}_0} \pair{dF(\rho),\hat\eta_0\rho} = \int_{\mathcal{B}_0} \eta_0 \cdot (\rho\times dF(\rho) )= \eta_0 \cdot (u_1+u_2). 
\end{align*}
Similarly, the bodies $\mathcal{B}_1$ and $\mathcal{B}_2$ are subject to the control moments $-R_1^TR_0u_1$ and $-R_2^T R_0 u_2$, respectively, with respect to their body fixed frames. The effect of the control moments is given by
\begin{align}
U = [u_1+u_2;\, 0_{3\times 1} ;\, -R_1^T R_0 u_1 ;\, -R_2^T R_0 u_2].\label{eqn:U}
\end{align}

Substituting \refeqn{U} into \refeqn{EL0}, and using the results of~\cite{LeeLeoPACC09}, the Euler-Lagrange equations for the connected rigid bodies immersed in a perfect fluid are given by
\begin{gather}
\begin{aligned}
&\mathbb{I}(R_0,R_1,R_2)
    \begin{bmatrix} \dot\Omega_0 \\ \ddot x \\ \dot\Omega_1 \\ \dot\Omega_2\end{bmatrix}\\
&+\begin{bmatrix}
    \Omega_0\times J\Omega_0 +\widehat{ R_0^T\dot x} M_0 R_0^T \dot x
    + \sum_{i=1}^2\hat d_{0i}R_0^T R_i W_i \\
    R_0(\hat\Omega_0M_0-M_0\hat\Omega_0)R_0^T\dot x +\sum_{i=1}^2 R_i W_i\\
    \Omega_1\times J_1\Omega_1 +V_1\times M_1V_1 -\hat d_{10}W_1\\
    \Omega_2\times J_2\Omega_2 +V_2\times M_2 V_2 -\hat d_{20}W_2
\end{bmatrix}=U,\end{aligned}\\
\dot R_0=R\hat\Omega_1,\quad \dot R_1=R_1\hat\Omega_1,\quad \dot R_2=R_2\hat\Omega_2,
\end{gather}
where
\begin{align}
V_i & = R_i^T \dot x - R_i^T R_0 \hat d_{0i} \Omega_0 + \hat d_{i0}\Omega_i,\\
W_i & =(\hat\Omega_iM_i-M_i\hat\Omega_i)(R_i^T\dot x - R_i^TR_0\hat d_{0i}\Omega_0)\nonumber\\
&\quad -M_iR_i^TR_0 \hat\Omega_0 \hat d_{0i}\Omega_0 + \hat\Omega_i M_i \hat d_{i0}\Omega_i
\end{align}
for $i\in\{1,2\}$.

\subsection{Symmetry}

Let the momentum of the system be $\mu=[p_0;p_x;p_1;p_2]\in\Re^{12}\simeq \g^*$. The Legendre transformation can be written as $\mu = \D_\xi L(g,\xi) = \mathbb{I}(R_0,R_1,R_2)\xi$. The total linear momentum of the connected rigid bodies and the ambient fluid is given by $p_x$, and the total angular momentum is given by $p_\Omega=\hat x p_x +\sum_{i=0}^2 R_i p_i \in\Re^3$.

The Lagrangian for the connected rigid bodies in a perfect fluid is left-invariant under rigid translation and rotation of the entire system. As a result, the configuration manifold can be reduced to a shape space $\G/\SE=\SO\times\SO$, and the total angular momentum and the total linear momentum are preserved. 

The control moments $u_1$ and $u_2$ change the relative attitudes between rigid bodies in the shape space $\SO\times\SO$. Since they respect the symmetry of the Lagrangian, the total linear momentum and the total angular momentum are also preserved in the controlled dynamics of the connected rigid bodies. 

\section{Computational Geometric Optimal Control of Connected Rigid Bodies in a Perfect Fluid}

Geometric numerical integration deals with numerical integrators that preserve geometric features of a dynamical system, such as invariants, symmetry, and reversibility~\cite{Hai.BK2006}. For numerical simulation of Hamiltonian systems that evolve on a Lie group, such as our system of connected rigid bodies in a fluid, it is critical to preserve both the symplectic property of Hamiltonian flows and the Lie group structure~\cite{Lee.CMaDA2007}. A geometric numerical integrator, referred to as a Lie group variational integrator, has been developed for a system of connected rigid bodies in a perfect fluid~\cite{LeeLeoPACC09}. It has desirable properties, such as preserving the symplecticity, momentum map, and Lie group structure of the system, thereby providing qualitatively correct numerical results for complex maneuvers over a long time period. 

The computational geometric optimal control approach utilizes the structure-preserving properties of geometric integrators in an optimization process~\cite{LeeLeoCISSIDRWB09,Lee.2008}. More explicitly, a discrete-time optimal control problem is formulated using Lie group variational integrators, and general optimization techniques, such as an indirect method based on optimality conditions or a direct method based on nonlinear programing, are applied. 

This method has substantial computational advantages. Since variational integrators compute the energy dissipation of a controlled system accurately~\cite{MarWesAN01}, more accurate optimal trajectories are obtained. Since there is no artificial numerical dissipation induced by the integration algorithms, this approach is numerically more robust, and the optimal control input can be computed efficiently. By representing the configuration of the rigid bodies directly on a Lie group, this method avoids singularities and complexities that appear in local coordinates. These properties are particularly useful for connected rigid bodies that are controlled indirectly by nonlinear coupling effects through large-angle rotational maneuvers.

\subsection{Lie Group Variational Integrator}

Here, we generalize the Lie group variational integrator for a system of connected rigid bodies in a perfect fluid developed in \cite{LeeLeoPACC09}, to include the effects of the internal control moments $u_1$ and $u_2$.

Let $h>0$ be a fixed integration step size, and let a subscript $k$ denote the value of a variable at the $k$-th time step. We define a discrete-time kinematic equation as follows. Define $f_k=(F_{0_k},\Delta x_k,F_{1_k},F_{2_k})\in\G$ for $\Delta x_k\in\Re^3$, $F_{0_k},F_{1_k},F_{2_k}\in\SO$ such that $g_{k+1}=g_k f_k$:
\begin{align}
    (R_{0_{k+1}},\,& x_{k+1},\,R_{1_{k+1}},\,R_{2_{k+1}})\nonumber\\
    & =(R_{0_{k}}F_{0_k},\,x_{k}+\Delta x_k,\,R_{1_{k}}F_{1_k},\,R_{2_{k}}F_{2_k}).
\end{align}
Therefore, $f_k$ represents the relative update between two integration steps. This ensures that the structure of the Lie group configuration manifold is numerically preserved.

\paragraph*{Discrete Lagrangian}

A discrete Lagrangian $L_d(g_k,f_k):\G\times\G\rightarrow\Re$ is an approximation of the Jacobi solution of the Hamilton--Jacobi equation, which is given by the integral of the Lagrangian along the exact solution of the Euler-Lagrange equations over a single time step. The discrete Lagrangian of connected rigid bodies is chosen as
\begin{align}
    L_{d_k} & = \frac{1}{2h}\Delta x_k^T R_{0_k}M_0 R_{0_k}^T\Delta x_k +\frac{1}{h}\tr{(I-F_{0_k})  J_{d_0}}\nonumber\\
    &\quad  +\sum_{i=1}^2 \Big( \frac{1}{2h}\Delta x_k^T
    R_{i_k}M_i R_{i_k}^T\Delta x_k + \frac{1}{h}\tr{ (I-F_{i_k}) J_{d_i}'}\nonumber\\
    &\quad + \frac{1}{2h} d_{0i}^T(F_{0_k}^T-I)R_{0_k}^T R_{i_k} M_i R_{i_k}^T R_{0_k} (F_{0_k}-I)d_{0i}\nonumber\\
    &\quad +\frac{1}{h}\Delta x_k^T R_{i_k} M_i R_{i_k}^T R_{0_k}(F_{0_k}-I)d_{0i}\nonumber\\
    &\quad -\frac{1}{h}\Delta x_k^T R_{i_k} M_i (F_{i_k}-I)d_{i0}\nonumber\\
    &\quad-\frac{1}{h}d_{0i}^T (F_{0_k}^T -I) R_{0_k}^T R_{i_k} M_i (F_{i_k}-I) d_{i0}\Big),\label{eqn:Ldk}
\end{align}
where nonstandard inertia matrices are defined as
\begin{gather}
    J_{d_0} = \frac{1}{2}\tr{J_0}I-J_0,\label{eqn:Jd0}\\
    J_{d_i}' = \frac{1}{2}\tr{J_i'}I-J_i',\quad
    J_i' = J_i - \hat d_{i0}M_i \hat d_{i0},\label{eqn:Jdi}
\end{gather}
for $i\in\{1,2\}$.

\paragraph*{Effects of Control Moments}

For a discrete Lagrangian defined on an arbitrary Lie group, the following discrete-time Euler-Lagrange equations, referred to as a Lie group variational integrator, were developed in~\cite{Lee.2008}.
\begin{gather}
\begin{aligned}
    \T^*_e & \L_{f_{k}}\cdot \D_{f_k}L_{d_k}-\Ad^*_{f_{k+1}^{-1}}\cdot(\T^*_e\L_{f_{k+1}}\cdot \D_{f_{k+1}}L_{d_{k+1}})\\
    &\quad+\T^*_e\L_{g_{k+1}}\cdot \D_{g_{k+1}} L_{d_{k+1}}+U_{d_k}^-+U_{d_{k-1}}^+=0,
\end{aligned}\label{eqn:DEL0}\\
    g_{k+1} = g_k f_k,\label{eqn:DEL1}
\end{gather}
where $\Ad^*:\G\times\g^*\rightarrow\g^*$ is the $\mathrm{co}$-$\mathrm{Adjoint}$ operator~\cite{Mar.BK1999}. The discrete generalized forces $U_{d_k}^+,U_{d_k}^-\in\g$ are chosen so as to approximate the virtual work of external forces and moments over a time step:
\begin{align*}
\int_{t_k}^{t_{k+1}} U(t)\cdot \eta \,dt \approx U_{d_k}^- \cdot \eta_k + U_{d_k}^+\cdot \eta_{k+1},
\end{align*}
where the variation of the configuration variable is represented by $\delta g =g \eta$ for $\eta\in\g$.

From \refeqn{U}, these are chosen as
\begin{gather}
U_{d_k}^- = \frac{h}{2} U_k,\qquad  
U_{d_k}^+ = \frac{h}{2} U_{k+1}.\label{eqn:Udk}
\end{gather}

\paragraph*{Discrete-time Euler-Lagrange Equations}

Substituting \refeqn{Ldk}, \refeqn{Udk} into \refeqn{DEL0}, \refeqn{DEL1}, the discrete-time Euler-Lagrange equations for the connected rigid bodies immersed in a perfect fluid are given by
\begin{gather}
\begin{aligned}
(J_{0_d}&F_{0_k}- F_{0_k}^TJ_{0_d})^\vee-(F_{0_{k+1}}J_{0_d}- J_{0_d}F_{0_{k+1}}^T)^\vee\\
&+(M_0 R_{0_{k+1}}^T \Delta x_{k+1})^\wedge R_{0_{k+1}}^T\Delta x_{k+1}\\
&+ \sum_{i=1}^2 \hat d_{0i} R_{0_{k+1}}^T (A_{i_k}- A_{i_{k+1}})=-h(u_{1_k}+u_{2_k}),
\end{aligned}\label{eqn:F0kp}\\
\begin{aligned}
(J'_{i_d}&F_{i_k}- F_{i_k}^TJ'_{i_d})^\vee
-(F_{i_{k+1}}J'_{i_d}- J'_{i_d}F_{i_{k+1}}^T)^\vee\\
&- \hat d_{i0}F_{i_k}^T  M_i R_{i_k}^TB_{i_k}+\widehat {F_{i_{k+1}} d_{i0}}  M_i R_{i_{k+1}}^TB_{i_{k+1}}\\
&+R_{i_{k+1}}^T \hat A_{i_{k+1}}B_{i_{k+1}}=hR_{i_k}^T R_{0_k} u_{i_k},
\end{aligned}\\
\begin{aligned}
R_{0_k}&M_0 R_{0_k}^T \Delta x_k + A_{1_k} + A_{2_k}\\
&-R_{0_{k+1}}M_0 R_{0_{k+1}}^T \Delta x_{k+1} - A_{1_{k+1}} - A_{2_{k+1}}=0,\\
\end{aligned}\label{eqn:delxkp}\\
R_{0_{k+1}}=R_{0_k} F_{0_k},\quad R_{i_{k+1}}=R_{i_k} F_{i_k}\label{eqn:R0kp}\\
x_{k+1} = x_k + \Delta x_k,\label{eqn:xkp}
\end{gather}
where inertia matrices are given by \refeqn{Jd0}, \refeqn{Jdi}, and $A_{i_k},B_{i_k}\in\Re^{3}$ are given by 
\begin{align}
A_{i_k} &= R_{i_k} M_i \parenth{R_{i_k}^T  B_{i_k} - (F_{i_k}-I)d_{i0}},\label{eqn:Aik}\\
B_{i_k} & = \Delta x_k + R_{0_k} (F_{0_k}-I)d_{0i}.\label{eqn:Bik}
\end{align}
For given $(g_0,f_0)\in\G\times\G$, $g_1\in\G$ is obtained by \refeqn{R0kp}--\refeqn{xkp}, and $f_1\in\G$ is obtained by solving \refeqn{F0kp}--\refeqn{delxkp}. This yields a discrete-time Lagrangian flow map $(g_0,f_0)\rightarrow(g_1,f_1)$, and this process is iterated.

This integrator preserves the symplectic form, momentum maps, and Lie group structure of a system of connected rigid bodies in a fluid exactly (up to machine precision), and the total energy conservation error is uniformly bounded over a long-time period when there is no control moment.

\subsection{Optimal Control Problem}

The objective of the optimal control problem of this paper is to transfer the connected rigid bodies from an initial configuration to a desired terminal condition during a fixed maneuver time $Nh$, while minimizing the control inputs:
\begin{gather}
\min_{u_k} \braces{\sum_{k=0}^{N} \frac{h}{2} (u_{1_k}\cdot u_{1_k} + u_{2_k}\cdot u_{2_k})}.\label{eqn:cost}
\end{gather}

As discussed before, the control moments $u_1$ and $u_2$ change the relative attitudes between rigid bodies in the shape space $\SO\times\SO$, and the total linear momentum and the total angular momentum are preserved in the controlled dynamics of the connected rigid bodies. 

A system of connected rigid bodies can translate and rotate as a consequence of the effects of geometric phase and momentum exchange between the rigid bodies and the ambient fluid. Geometric phase refers to a translation along the symmetry direction achieved by closed trajectories in the shape space~\cite{Mar92}. By controlling the relative attitudes, a system of connected rigid bodies can translate and rotate without external forces or moments. Based on geometric phase effects, the controllability and motion planning algorithms for articulated multibody systems with reaction wheels have been studied in~\cite{ReyMcCAJGCD92,RuiKolITAC00}, and several optimal control problems for multibody systems have been studied in~\cite{LeeLeoPACC07,LeeLeoCISSIDRWB09}. In addition to geometric phase, connected rigid bodies can steer themselves by exchanging linear momentum and angular momentum with the fluid. 

\subsection{Computational Approach}

We apply a direct optimal control approach. The control inputs are parameterized by several points that are uniformly distributed over the maneuver time, and control inputs between these points are approximated using cubic spline interpolation. For given control input parameters, the value of the cost is given by \refeqn{cost}, and the terminal conditions are obtained by the discrete-time equations of motion  \refeqn{F0kp}-\refeqn{xkp}. The control input parameters are optimized using constrained nonlinear parameter optimization to satisfy the terminal boundary conditions while minimizing the cost.

This approach is computationally efficient when compared to the usual collocation methods, where the continuous-time equations of motion are imposed as constraints at a set of collocation points. Using the proposed discrete-time optimal control approach, optimal control inputs can be obtained by using a large step size, thereby resulting in efficient total computations. Since the computed optimal trajectories do not have numerical dissipation caused by conventional numerical integration schemes, they are numerically more robust. Furthermore, the corresponding gradient information of the cost and the terminal constraints with respect to the control parameters is accurately computed, which improves the convergence properties of the numerical optimization procedure.

\section{Numerical Example}\label{sec:NE}

The properties of connected rigid bodies are chosen as follows. The principal axes of each ellipsoid are given by
\begin{align*}
\text{Body 0: } & l_1=8,\quad l_2 =1.5,\quad l_3=2\;(\mathrm{m}),\\
\text{Body 1,2: } & l_1=5,\quad l_2 =0.8,\quad l_3=1.5\;(\mathrm{m}).
\end{align*}
We assume the density of the fluid is $\rho=1\mathrm{kg/m^3}$. The corresponding inertia matrices are given by
\begin{gather*}
    M_0= \mathrm{diag}[1.0659,\,2.1696,\,1.6641],\;(\mathrm{kg})\\
    M_1=M_2=\mathrm{diag}[0.2664,\,    0.6551,\,    0.3677]\;(\mathrm{kg}),\\
    J_0= \mathrm{diag}[1.3480,\,   20.1500,\,   25.3276]\;(\mathrm{kgm^2}),\\
    J_1=J_2=\mathrm{diag}[0.1961,\,    1.7889,\,    2.9210]\;(\mathrm{kgm^2}).
\end{gather*}
The location of the ball joints with respect to the center of mass of each body are chosen as
\begin{gather*}
d_{01}=-d_{02}=[8.8,\,0,\,0],\quad d_{10}=-d_{20}=[5.5,\,0,\,0]\;(\mathrm{m}).
\end{gather*}
The initial configuration is as follows:
\begin{gather*}
    R_{0_0}= R_{1_0}=I,\quad
    R_{2_0}=\exp\parenth{\frac{\pi}{4}\hat e_3},\quad
    x_0=0_{3\times 1}\;(\mathrm{m}).
\end{gather*}
The initial velocities are set to zero, i.e $\dot x_0,\Omega_0,\Omega_1,\Omega_2=0_{3\times 1}$.

We consider the following two rest-to-rest maneuvers. 
\begin{itemize}
\item[(i)] Forward translation along the $e_1$ axis
\begin{gather*}
e_1\cdot x_N =2,\qquad R_{i_0}^TR_{i_N}=I\quad \text{for $i\in\{0,1,2\}$}.
\end{gather*}
\item[(ii)] $180^\circ$ rotation about the $e_1$ axis
\begin{gather*}
R_{i_N}=\exp (\pi\hat e_1) R_{i_0} \quad \text{for $i\in\{0,1,2\}$}.
\end{gather*}
\end{itemize}
In every case, the terminal velocities are set to zero, and the unspecified parts of the terminal position are free. The maneuver time is $T=1$ second, and the time step is $h=0.001$ second. The first case is a planar motion, where the connected rigid bodies move in the $e_1e_2$ plane, and the control moments $u_1,u_2$ act along the $e_3$ axis. The second case is a three-dimensional maneuver. 

\begin{figure}
\subfigure[Snapshots at each $0.7$ second]{
\includegraphics[width=\columnwidth]{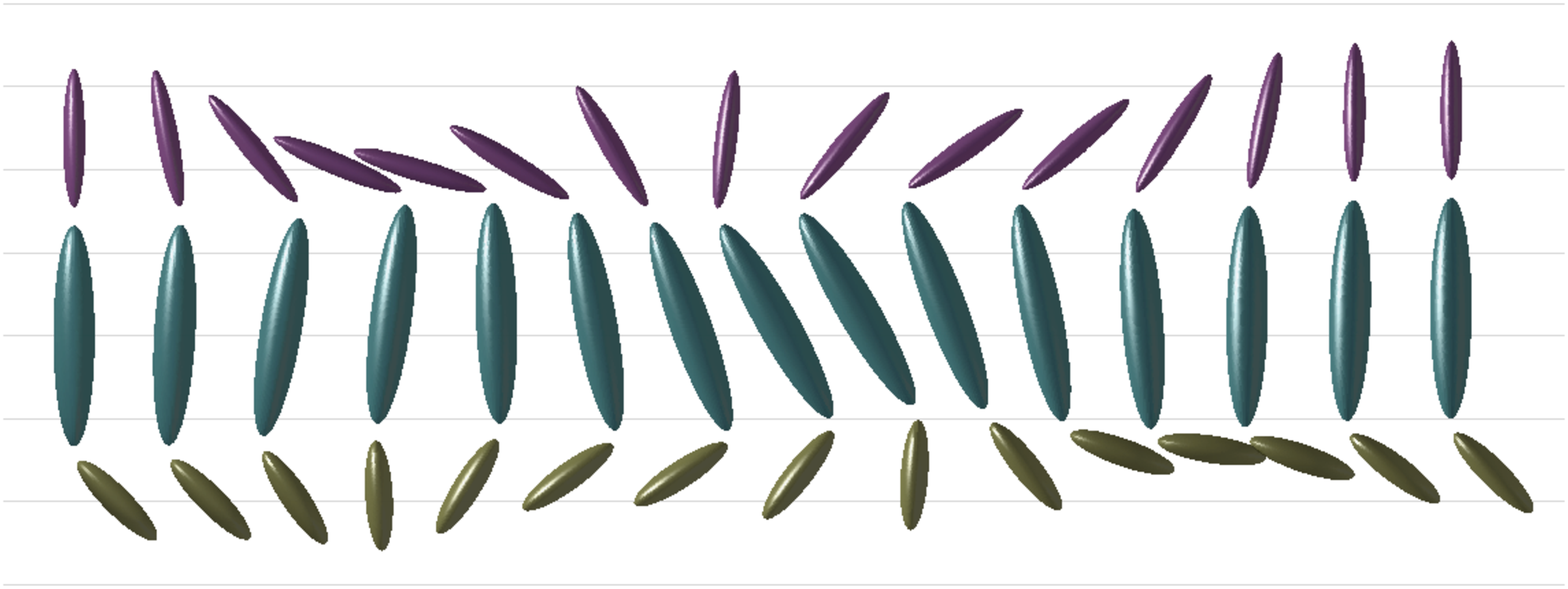}}
\centerline{
\subfigure[Velocity of the central body (first and second component ($\mathrm{m/s}$))]
	{\includegraphics[width=0.49\columnwidth]{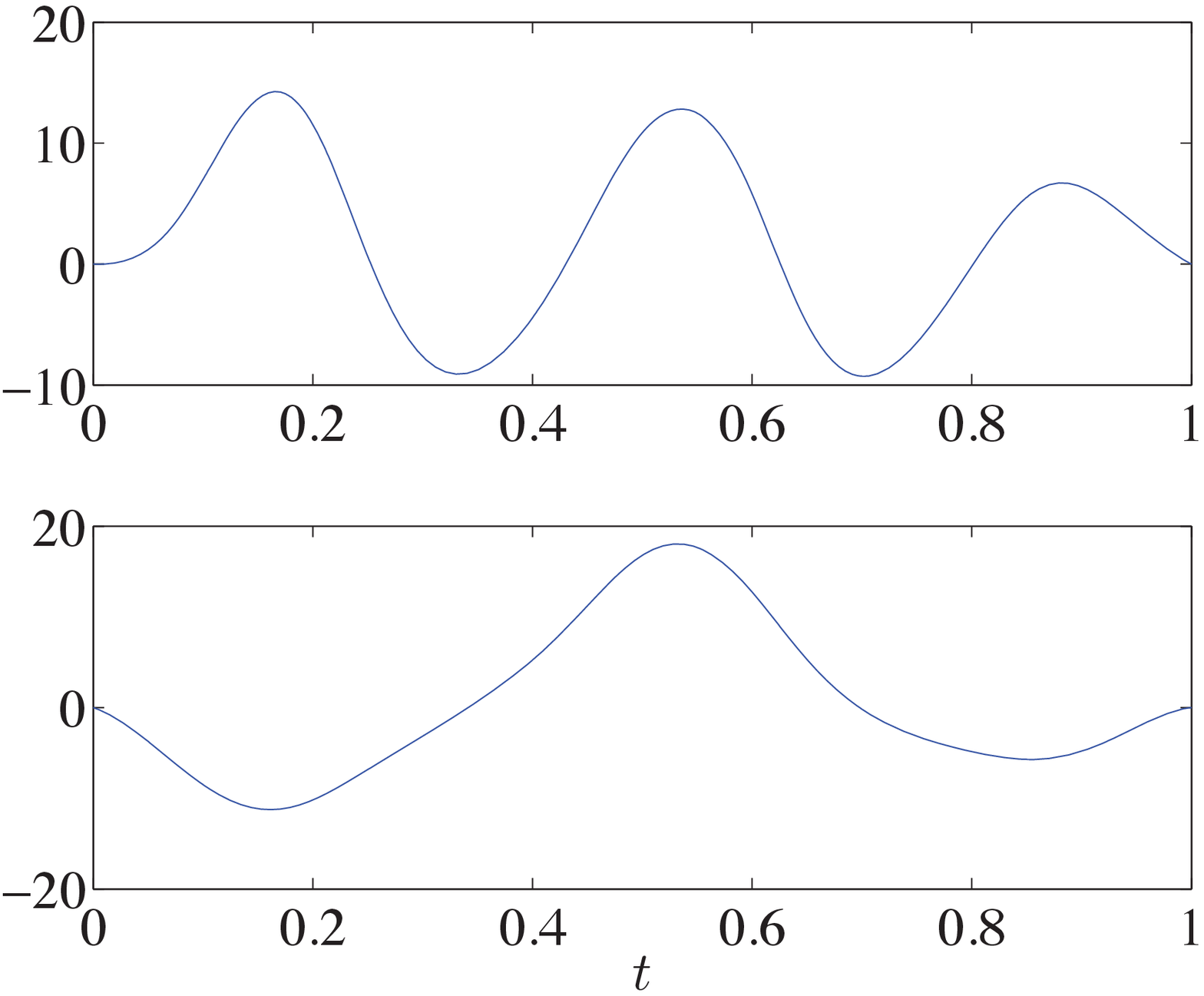}}
\hspace*{0.02\columnwidth}
\subfigure[Angular velocity (third components, red:$\Omega_1$, blue:$\Omega_0$, green:$\Omega_2$, ($\mathrm{rad/s}$))]
	{\includegraphics[width=0.49\columnwidth]{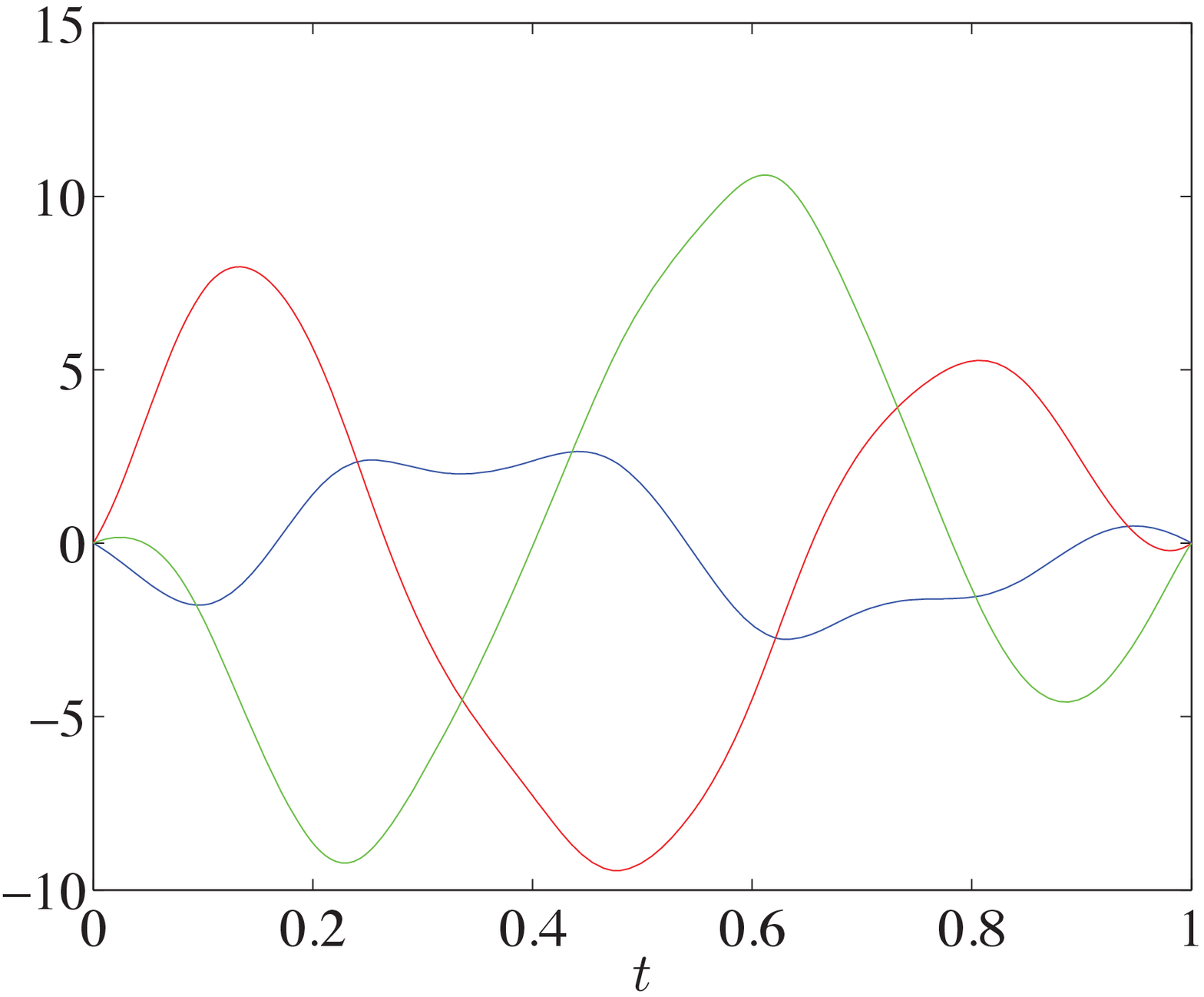}}
}
\centerline{
\subfigure[Linear momentum along the $e_1$ axis (dashed:rigid bodies, dotted:fluid, solid:total ($\mathrm{kgm/s}$))]
	{\includegraphics[width=0.49\columnwidth]{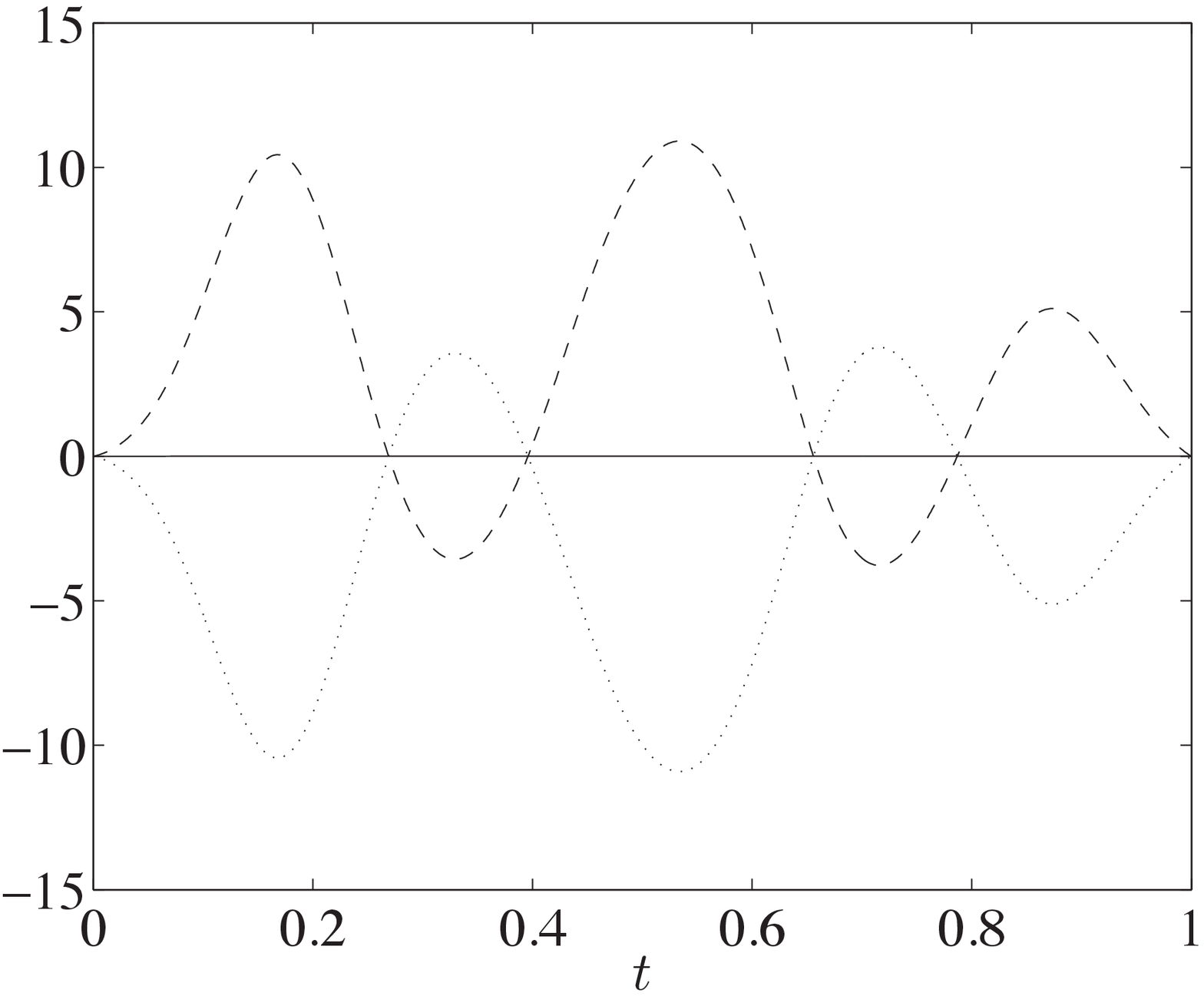}\label{fig:i_px}}
\hspace*{0.02\columnwidth}
\subfigure[Control moment (third components, red:$u_1$, green:$u_2$, ($\mathrm{kNm}$))]
	{\includegraphics[width=0.49\columnwidth]{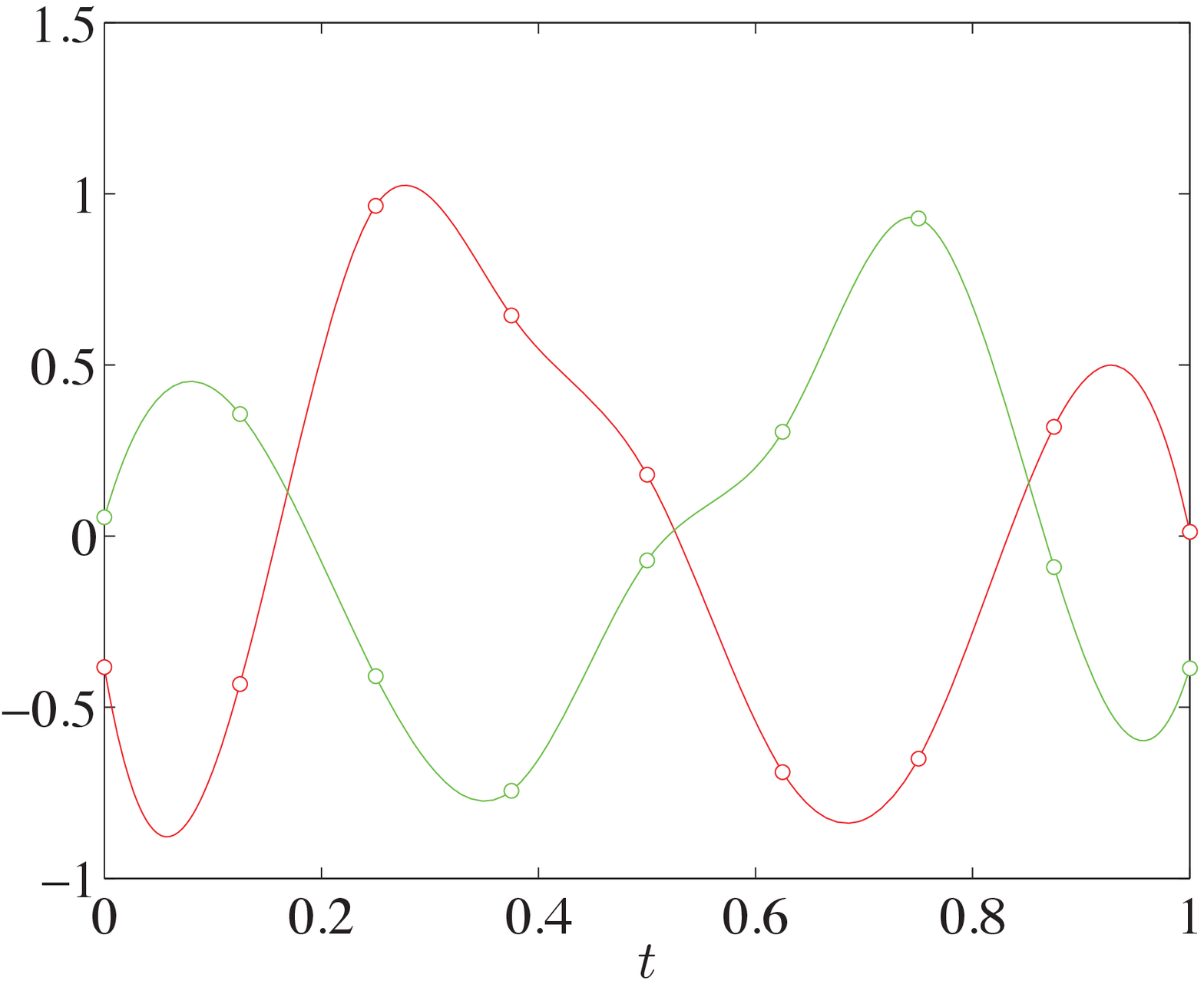}}
}
\caption{Optimal forward translation}\label{fig:i}
\end{figure}

\reffig{i} shows the optimal forward translation maneuver (an animation illustrating this maneuver is available at {\small\url{http://my.fit.edu/~taeyoung}}). The linear momentum exchange between rigid bodies and fluid along the $e_1$ axis is illustrated at \reffig{i_px}. While the zero value of the total linear momentum is preserved throughout the maneuver, the average value of the linear momentum of the rigid bodies is positive. Therefore, the system of connected rigid bodies moves forward.

\begin{figure}
\subfigure{\includegraphics[width=\columnwidth]{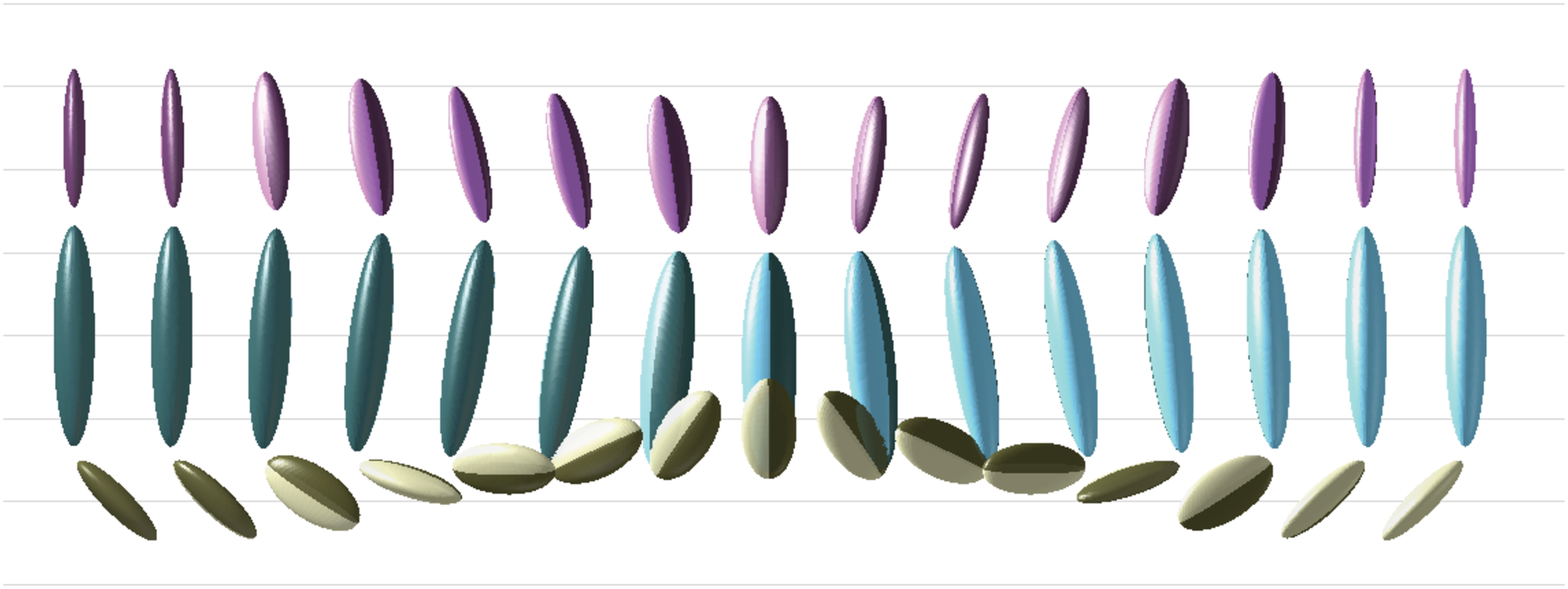}}
\setcounter{subfigure}{0}
\subfigure[Snapshots at each $0.7$ second (top view and back view)]
	{\includegraphics[width=\columnwidth]{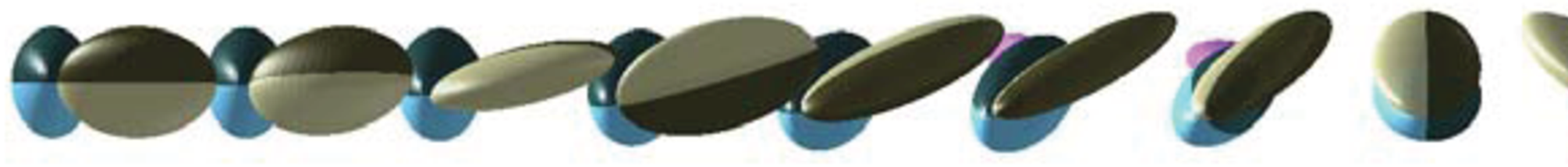}}
\centerline{
\subfigure[Velocity of the central body (each component ($\mathrm{m/s}$))]
	{\includegraphics[width=0.49\columnwidth]{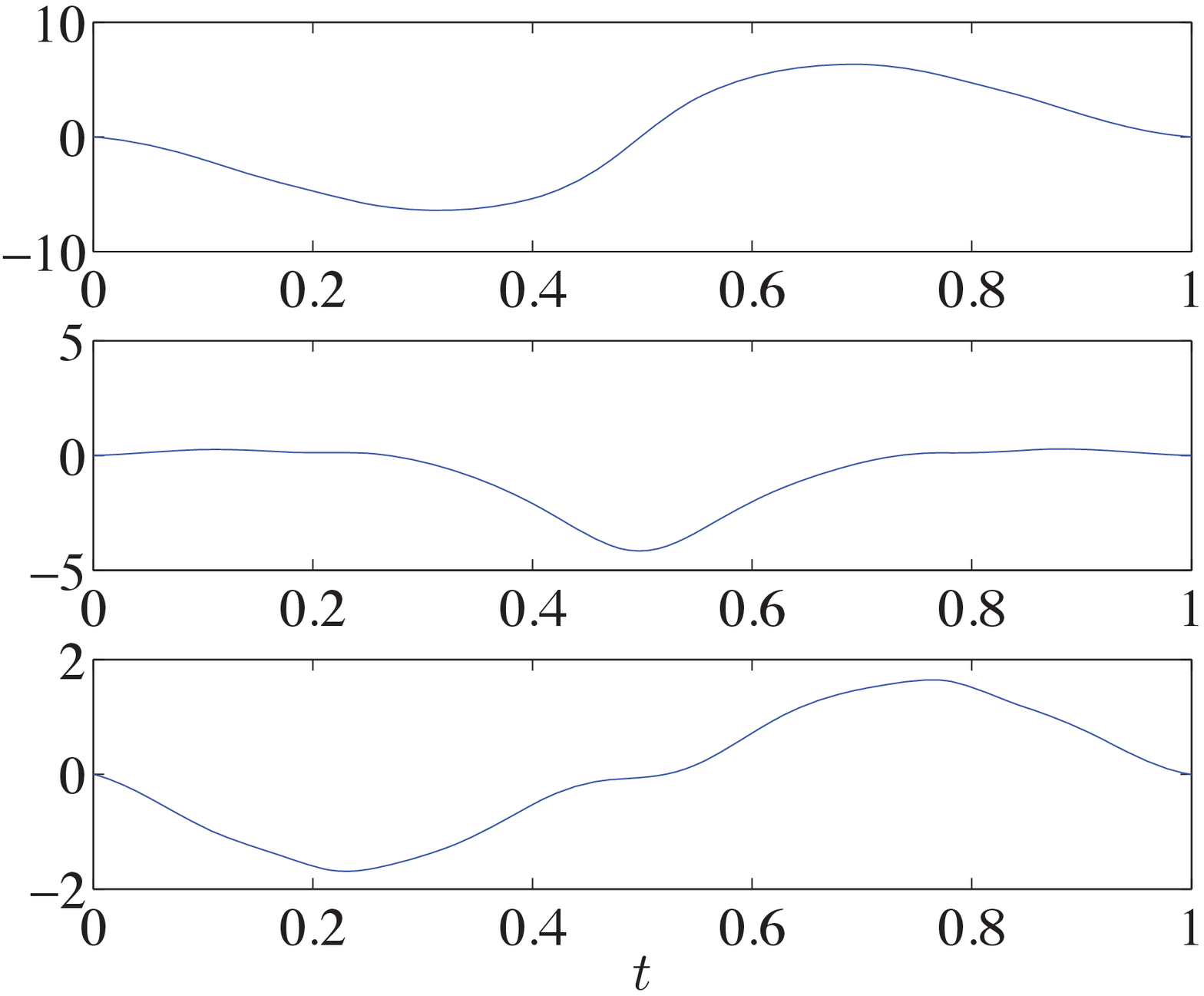}}
\hspace*{0.02\columnwidth}
\subfigure[Angular velocity (each component, red:$\Omega_1$, blue:$\Omega_0$, green:$\Omega_2$, ($\mathrm{rad/s}$))]
	{\includegraphics[width=0.49\columnwidth]{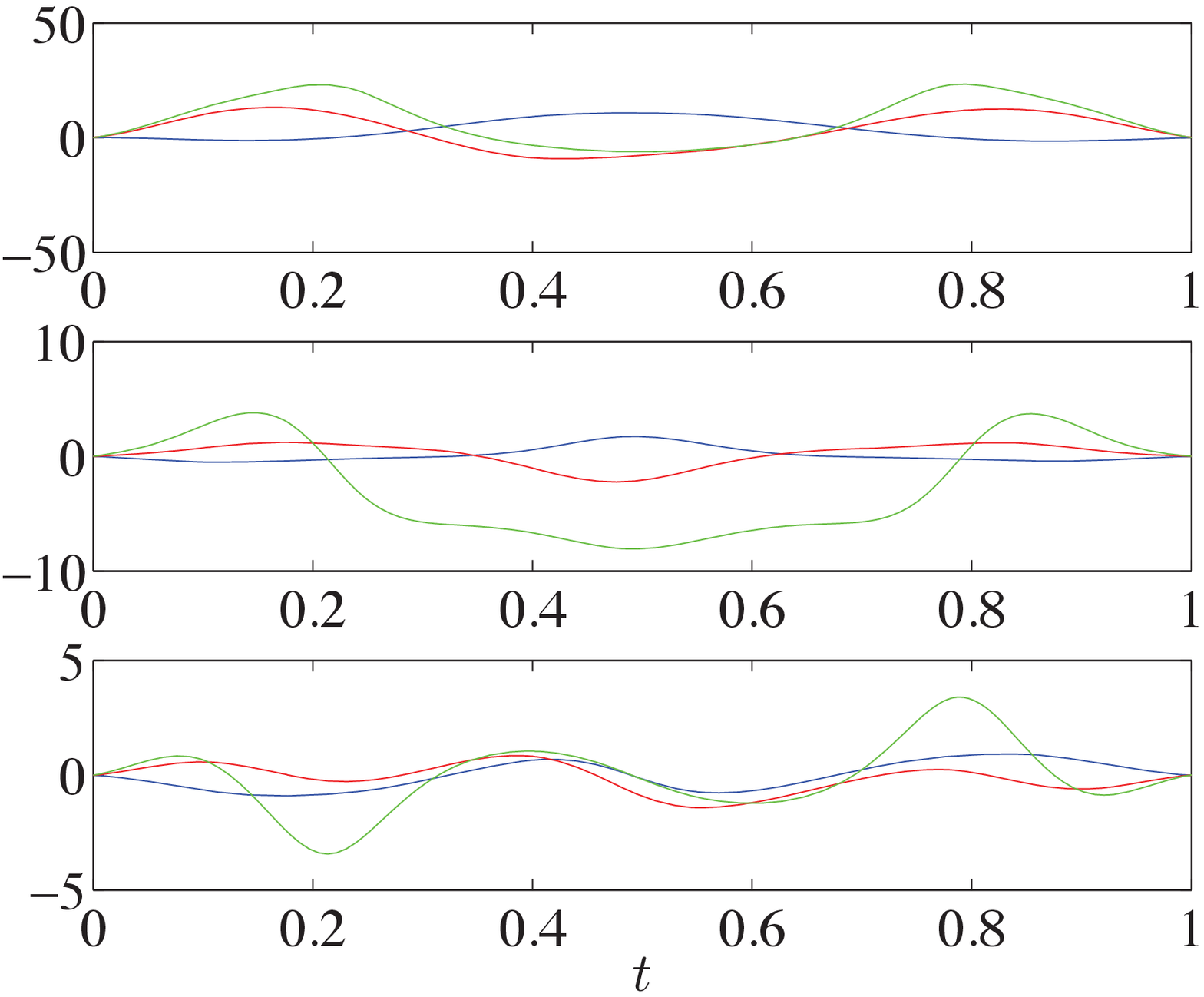}}
}
\centerline{
\subfigure[Angular momentum about the $e_1$ axis (dahsed:rigid bodies, dotted:fluid, solid:total ($\mathrm{kgm/s}$))]
	{\includegraphics[width=0.48\columnwidth]{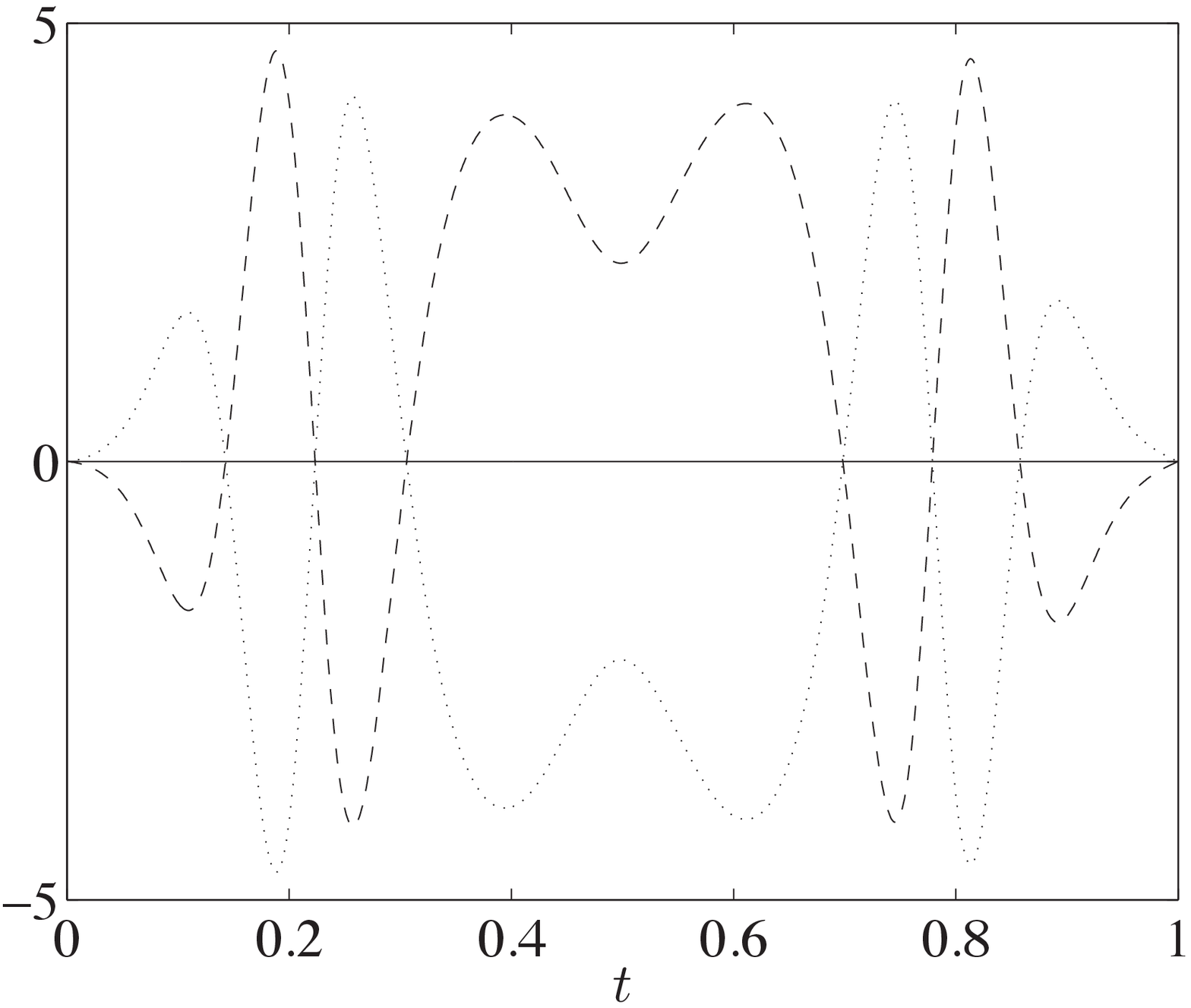}\label{fig:ii_pw}}
\hspace*{0.02\columnwidth}
\subfigure[Control moment (each component, red:$u_1$, green:$u_2$, ($\mathrm{kNm}$))]
	{\includegraphics[width=0.50\columnwidth]{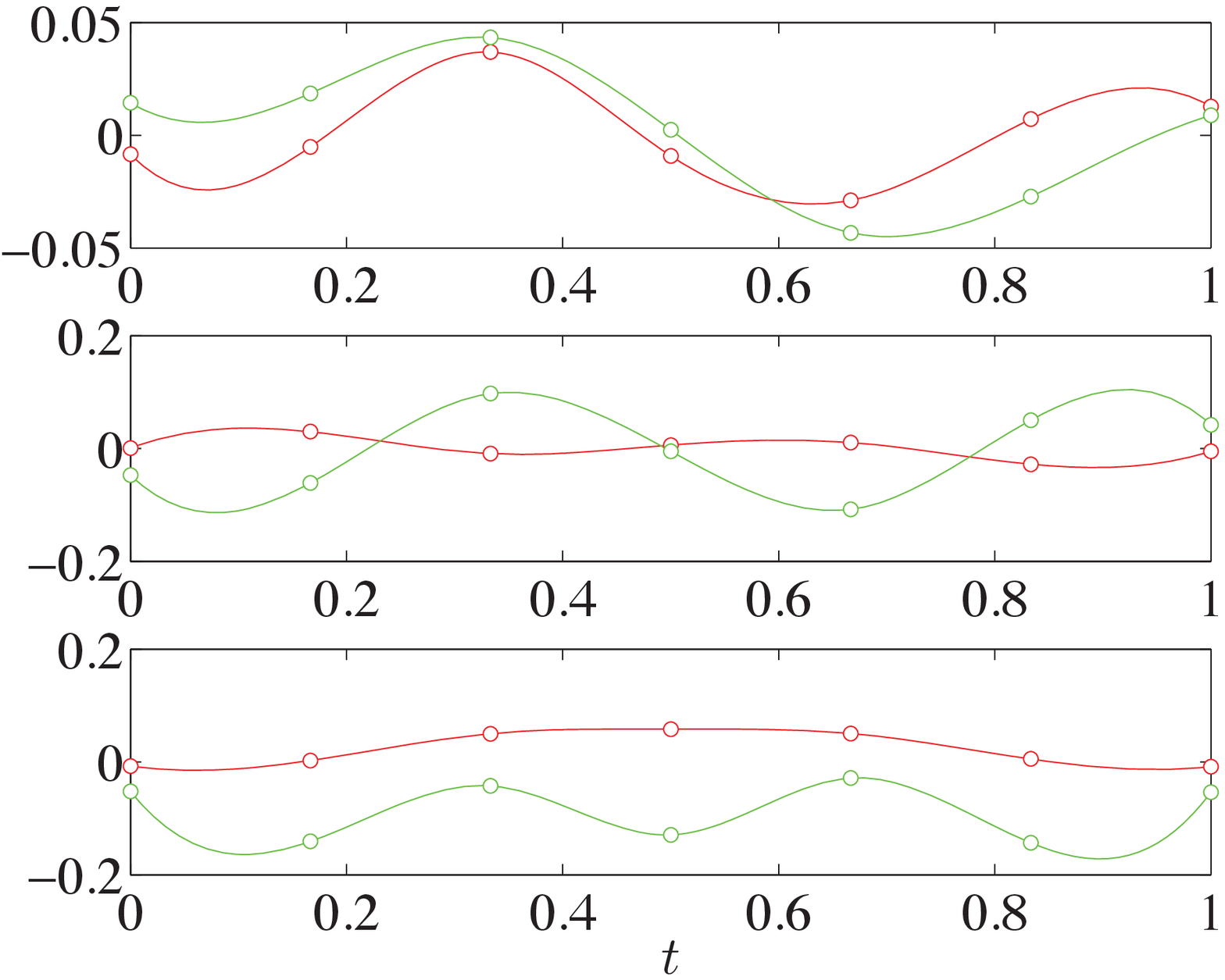}}
}
\caption{Optimal $180^\circ$ rotation about the $e_1$ axis}\label{fig:ii}
\end{figure}

\reffig{ii} illustrates the optimal rotation maneuver. It is interesting to observe that the optimal velocity, angular velocity, and control inputs are almost symmetric about the mid-maneuver time $t=0.5$. Similar to the optimal translation maneuver, at \reffig{ii_pw}, the total angular momentum is preserved, but the average angular momentum of the rigid bodies about the rotation axis is positive. In addition to the angular momentum exchange between the rigid bodies and the fluid, this rotational maneuver is achieved by the geometric phase effect for coupled rigid bodies~\cite{MonDCMSFRP91}.

\section{Conclusion}
For both cases, the proposed computational geometric optimal control approach successfully finds nontrivial optimal maneuvers of connected rigid bodies in a perfect fluid. These optimal maneuvers are based on the nonlinear coupling between the shape space and the group motions, and they involve large-angle rotations of connected rigid bodies.

The computational framework described in this paper extends the class of computationally tractable shape-based fish locomotion models to include three-dimensional articulated multibody systems, thereby allowing for the study of optimal shape maneuvers for fish movement patterns beyond that of carangiform locomotion.

\bibliography{ACC10}

\begin{thebibliography}{10}
\providecommand{\url}[1]{#1}
\csname url@rmstyle\endcsname
\providecommand{\newblock}{\relax}
\providecommand{\bibinfo}[2]{#2}
\providecommand\BIBentrySTDinterwordspacing{\spaceskip=0pt\relax}
\providecommand\BIBentryALTinterwordstretchfactor{4}
\providecommand\BIBentryALTinterwordspacing{\spaceskip=\fontdimen2\font plus
\BIBentryALTinterwordstretchfactor\fontdimen3\font minus
  \fontdimen4\font\relax}
\providecommand\BIBforeignlanguage[2]{{%
\expandafter\ifx\csname l@#1\endcsname\relax
\typeout{** WARNING: IEEEtran.bst: No hyphenation pattern has been}%
\typeout{** loaded for the language `#1'. Using the pattern for}%
\typeout{** the default language instead.}%
\else
\language=\csname l@#1\endcsname
\fi
#2}}

\bibitem{Sfa.IJoOE1999}
M.~Sfakiotakis, D.~Lane, and J.~Davies, ``Review of fish swimming modes for
  aquatic locomotion,'' \emph{IEEE Journal of Oceanic Engineering}, vol.~24,
  no.~2, pp. 237--252, 1999.

\bibitem{Bar.1996}
D.~Barrett, ``Propulsive efficiency of a flexible hull underwater vehicle,''
  Ph.D. dissertation, Massachusetts Institute of Technology, 1996.

\bibitem{Kan.JoNS2005}
E.~Kanso, J.~Marsden, C.~Rowley, and J.~Melli-Huber, ``Locomotion of
  articulated bodies in a perfect fluid,'' \emph{Journal of Nonlinear Science},
  vol.~15, pp. 255--289, 2005.

\bibitem{LeeLeoPACC09}
\BIBentryALTinterwordspacing
T.~Lee, M.~Leok, and N.~H. McClamroch, ``Dynamics of connected rigid bodies in
  a perfect fluid,'' in \emph{Proceedings of the American Control Conference},
  2009, pp. 408--413. [Online]. Available: \url{http://arxiv.org/abs/0809.1488}
\BIBentrySTDinterwordspacing

\bibitem{Lee.2008}
T.~Lee, ``Computational geometric mechanics and control of rigid bodies,''
  Ph.D. dissertation, University of Michigan, 2008.

\bibitem{Lee.CMaDA2007}
T.~Lee, M.~Leok, and N.~H. McClamroch, ``Lie group variational integrators for
  the full body problem in orbital mechanics,'' \emph{Celestial Mechanics and
  Dynamical Astronomy}, vol.~98, no.~2, pp. 121--144, June 2007.

\bibitem{Mar92}
J.~Marsden, \emph{Lectures on Mechanics}, ser. London Mathematical Society
  Lecture Note Series 174.\hskip 1em plus 0.5em minus 0.4em\relax Cambridge
  University Press, 1992.

\bibitem{LeeLeoCISSIDRWB09}
\BIBentryALTinterwordspacing
T.~Lee, M.~Leok, and N.~H. McClamroch, ``Computational geometric optimal
  control of rigid bodies,'' \emph{Communications in Information and Systems,
  special issue dedicated to R. W. Brockett}, vol.~8, no.~4, pp. 445--472,
  2008. [Online]. Available: \url{http://arxiv.org/abs/0805.0639}
\BIBentrySTDinterwordspacing

\bibitem{Kan.PotICoDaC2005}
E.~Kanso and J.~Marsden, ``Optimal motion of an articulated body in a perfect
  fluid,'' in \emph{Proceedings of the IEEE Conference on Decision and
  Control}, 2005, pp. 2511--2516.

\bibitem{Ros.PotACC2006}
S.~Ross, ``Optimal flapping strokes for self-propulsion in a perfect fluid,''
  in \emph{Proceedings of the American Control Conference}, 2006, pp.
  4118--4122.

\bibitem{Hol.PD1998}
P.~Holmes, J.~Jenkins, and N.~Leonard, ``Dynamics of the {K}irchhoff equations
  {I}: Coincident centers of gravity and bouyancy,'' \emph{Physica D}, vol.
  118, pp. 311--342, 1998.

\bibitem{Mar.BK1999}
J.~Marsden and T.~Ratiu, \emph{Introduction to Mechanics and Symmetry},
  2nd~ed., ser. Texts in Applied Mathematics.\hskip 1em plus 0.5em minus
  0.4em\relax Springer-Verlag, 1999, vol.~17.

\bibitem{Hai.BK2006}
E.~Hairer, C.~Lubich, and G.~Wanner, \emph{Geometric Numerical Integration},
  2nd~ed., ser. Springer Series in Computational Mathematics.\hskip 1em plus
  0.5em minus 0.4em\relax Springer-Verlag, 2006, vol.~31.

\bibitem{MarWesAN01}
J.~Marsden and M.~West, ``Discrete mechanics and variational integrators,'' in
  \emph{Acta Numerica}.\hskip 1em plus 0.5em minus 0.4em\relax Cambridge
  University Press, 2001, vol.~10, pp. 317--514.

\bibitem{ReyMcCAJGCD92}
M.~Reyhanoglu and N.~H. McClamroch, ``Reorientation maneuvers of planar
  multibody systems in space using internal controls,'' \emph{AIAA Journal of
  Guidance, Control and Dynamics}, vol.~15, no.~6, pp. 1475--1480, 1992.

\bibitem{RuiKolITAC00}
C.~Rui, I.~Kolmanovsky, and N.~H. McClamroch, ``Nonlinear attitude and shape
  control of spacecraft with articulated appendages and reaction wheels,''
  \emph{IEEE Transactions on Automatic Control}, vol.~45, no.~8, pp.
  1455--1469, 2000.

\bibitem{LeeLeoPACC07}
\BIBentryALTinterwordspacing
T.~Lee, M.~Leok, and N.~H. McClamroch, ``Optimal attitude control for a rigid
  body with symmetry,'' in \emph{Proceedings of the American Control
  Conference}, 2007, pp. 1073--1078. [Online]. Available:
  \url{http://arxiv.org/abs/math.OC/06009482}
\BIBentrySTDinterwordspacing

\bibitem{MonDCMSFRP91}
R.~Montgomery, ``Guage theory of the falling cat,'' in \emph{Dynamics and
  control of mechanical systems: the falling cat and related problems}, ser.
  Fields Institute Communications, M.~Enos, Ed.\hskip 1em plus 0.5em minus
  0.4em\relax American Mathematical Society, 1991, pp. 193--218.

\end{thebibliography}
\bibliographystyle{IEEEtran}

\end{document}